\documentclass[12pt]{article}
\usepackage{amsmath}
\usepackage{amssymb}
\usepackage{amsthm}
\usepackage{amsfonts}
\usepackage{amscd}
\usepackage[mathscr]{eucal}
\newcommand{\Z} {{\mathbb  Z}}
\newcommand{\Q}{{\mathbb  Q}}

\newcommand{\R} {{\mathbb R}}

\begin{document}
\parindent  25pt
\baselineskip  10mm
\textwidth  15cm    \textheight  23cm \evensidemargin -0.06cm
\oddsidemargin -0.01cm

\title{{Lattice-ordered matrix algebras over real GCD-domains }}
\author{\mbox{} {Fei Li }  \\
{\small School of Mathematical Sciences,
Institute of Mathematics and Interdisciplinary Science,}\\
{\small Capital Normal University, Beijing 100048, P.R.China,} \\
{\small Current address: School of Statistics and Applied Mathematics,} \\
{\small Anhui University of Finance and Economics, bengbu 233030,
Anhui Province, P.R.China } \\ \\
{Xianlong Bai } \\
{\small School of Mathematical Sciences,
Institute of Mathematics and Interdisciplinary Science, } \\
{\small Capital Normal University, Beijing 100048, P.R.China } \\ \\
{Derong Qiu }
\thanks{ \quad Corresponding author, E-mail:
derong@mail.cnu.edu.cn } \\
{\small School of Mathematical Sciences,
Institute of Mathematics and Interdisciplinary Science, } \\
{\small Capital Normal University, Beijing 100048, P.R.China } }

\date{}
\maketitle
\parindent  24pt
\baselineskip  10mm
\parskip  0pt

\par     \vskip  0.4 cm

{\bf Abstract} \ Let $ R \subset \R $ be a GCD-domain.
In this paper, Weinberg's conjecture on the $ n \times n $
matrix algebra $ M_{n}(R) \ (n \geq 2) $ is proved. Moreover, all the lattice orders
(up to isomorphisms) on a full $ 2 \times 2 $ matrix algebra over $ R $ are obtained.

\par   \vskip 0.2cm


\par  \vskip 0.3cm

\hspace{-0.6cm}{\bf 1. \ Introduction}

\par \vskip 0.2 cm

Let $ A $ be a commutative $ l-$ring and $ M_{n}(A) $ be the $ n \times n $
matrix ring over $ A. $ Then $ M_{n}(A) $
becomes an $ l-$algebra over $ A $ with the $ usual \ lattice \ order $
associated to the positive cone $ P = M_{n}(A^{+}), $ denoted by
$ (M_{n}(A), M_{n}(A^{+})) $ (see [MR]).
In 1966, Weinberg in [W] conjectured that, for any integer $ n \geq 2 $ and
$ A = \Q $ (the field of rational numbers),
if $ M_{n}(\Q) $ is an $ l-$algebra in which the identity matrix is positive,
then \ $ (M_{n}(\Q), P) \cong (M_{n}(\Q), M_{n}(\Q^{+})).
$ Also in [W] he proved that this conjecture is true for the case $ n = 2 $
and determined all but one of the lattice orders of $ M_{2}(\Q). $ Later,
S. A. Steinberg studied Weinberg's conjecture over totally ordered
fields (see [St1]). In 2002, J. Ma and P. Wojciechowski proved Weinberg's
conjecture over a totally ordered subfield of the real number field $ \R $
(see[MW]). In 2007, J. Ma and R. H. Redfield proved Weinberg's
conjecture over the ring of integers $ \Z $ (see[MR]).

Let $ R \subset \R $ be a GCD-domain.
In this paper, Weinberg's conjecture on the $ n \times n $
matrix ring $ M_{n}(R) \ (n \geq 2) $ is proved (see Theorem 2.4 in the following)
using the procedures developed in [MR]. Moreover, all the lattice orders
(up to isomorphisms) on a full $ 2 \times 2 $ matrix
algebra over $ R $ are obtained (see Theorem 2.5 in the following).
\par  \vskip 0.3cm

\hspace{-0.6cm}{\bf 2. \ Lattice-ordered matrix algebras $ M_{n}(R) $ }

\par \vskip 0.2 cm

Let $ R \subset \R $ be a GCD-domain, i.e., a domain in which any two
non-zero elements have a greatest common divisor (see [K]),
and $ K $ be the field of fractions of $ R. $ Let $ M_{n}(R) $ and $ M_{n}(K)$
be the $ n \times n \ (n \geq 2) $ matrix ring over $ R $ and $ K $ respectively.
Now for the $ l-$algebra $ (M_{n}(R), P) $ with a
positive cone $ P $ (see [Bi], [St2]), we denote
$ \overline{P} = \{A\in M_{n}(K) \mid kA \in P \ \text{for some}
\ 0 < k \in R \}. $ Throughout this paper, $ P $
denotes a general positive cone on $ M_{n}(R); \ \overline{P} $ denotes the
positive cone on $ M_{n}(K) $ extended from $ P ; \ P_{A} = A M_{n}(R^{+}) $
denotes a positive cone on $ M_{n}(R) $ with $ A \in M_{n}(R^{+}); $
and $ P_{D} = D M_{n}(K^{+}) $ denotes a positive cone on $ M_{n}(K) $ with
$ D \in M_{n}(K^{+}) $ (see[MR], [St1]). We also denote the group of units of $ R $
by $ R^{\times } = \{r \in R : \ r t = 1 \ \text{for some} \ t \in R \}. $
\par \vskip 0.2 cm

{\bf Proposition 2.1.} \ If $ (M_{n}(R), P) $ is an $ l-$algebra
over $ R, $ then, as an $ l-$module over $ R, $ it has a $ vl-$ basis with $ n^{2} $
elements $ \{B_{ij} : \ 1 \leq i, j \leq n \}, $ and there exist two
non-singular matrices
$ H \in M_{n}(K), D = (d_{ij}) \in M_{n}(K^{+}), $ and a matrix
$ C = (q_{ij}) \in M_{n}(K^{+} \setminus \{\ 0 \}), $ such that the following
statements hold: \\
(1) \ $ B_{ij} = q_{ij} H D E_{ij} H^{-1} \ (1 \leq i, j \leq n); $ \\
(2) \ $ d_{jr} q_{ij} q_{rs} q_{is}^{-1} \in R^{+} \
(1 \leq i, j, r, s \leq n); $ \\
(3) \ $ B_{ij} B_{rs} = d_{jr} q_{ij} q_{rs} q_{is}^{-1} B_{is} \
(1 \leq i, j, r, s \leq n); $ \\
(4) \ $ (\prod_{1 \leq i, j \leq n} q_{ij} ) (det(D))^{n} \in R^{\times }. $ \\
In particular, if the identity matrix $ I \in P, $ then one can take
the above $ D = I, $ and then all $ q_{ij} q_{js} q_{is}^{-1} $ must be
positive units in $ R \ (1 \leq i, j, s \leq n). $
\par \vskip 0.1 cm
The proof of this proposition can be obtained from the proofs of Propositions 2.2, 2.3 and
Theorem 3.1 in [MR] by replacing the integers $ \Z $ by the real GCD-domain $ R. $
\par \vskip 0.2 cm

{\bf Theorem 2.2.} \ Let $ (M_{n}(R), P) $ be an
$ l-$algebra with a $ vl-$basis $ \{B_{ij} : \ 1 \leq i, j \leq n \} $
over $ R. $ If the identity matrix $ I \in P, $ then
$ (M_{n}(R), P) \cong (M_{n}(R), M_{n}(R^{+})) $
if and only if the system of equations $ x_{ij} x_{js} x_{is}^{-1} =
q_{ij} q_{js} q_{is}^{-1} \ (1 \leq i, j, s \leq n) $ with variables
$ x_{i^{\prime } j^{\prime }} \ (1 \leq i^{\prime }, j^{\prime } \leq n) $
has positive solutions in $ R^{\times }. $ Here
$ q_{ij} \ (1 \leq i, j \leq n ) $ are as in Proposition 2.1
above for the case $ D = I. $
\par \vskip 0.1 cm
{\bf Proof.} \ Since $ I \in P, $ by taking $ D = I $ in Prop. 2.1 above,
we have $ B_{ij} B_{js} = q_{i j} q_{j s} q_{i s}^{-1} B_{is}. $ Assume
$ (M_{n}(R), P) \cong (M_{n}(R), M_{n}(R^{+})), $ we denote such an isomorphism
by $ \psi . $ Then $ \{ \psi (B_{ij}) : \ 1 \leq i, j \leq n \} $ is a $ vl-$basis of
$ (M_{n}(R), M_{n}(R^{+})) $ over $ R. $ Since
$ \{E_{i j} : \ 1 \leq i, j \leq n \} $ is also a $ vl-$basis of
$ (M_{n}(R), M_{n}(R^{+})) $ over $ R, $ we may assume that $ \psi (B_{ij})
= \mu _{i j} E_{u(i,j)v(i,j)}, $ where $ \mu _{i j} \in R^{\times} \cap R^{+}. $
So $ \psi (B_{ij}) \psi (B_{js}) = q_{i j} q_{j s} q_{i s}^{-1}
\psi (B_{is}). $ Hence
$ \mu _{i j} \mu _{j s} E_{u(i,j)v(i,j)} E_{u(j,s)v(j,s)} =
q_{i j} q_{j s} q_{i s}^{-1} \mu _{is} E_{u(i,s)v(i,s)}. $ If
$ i = j = s, $ then we have $ u(i,i) = v(i,i) $ and $ \mu _{ii} = q_{ii}. $
So there exist an $ \sigma \in S_{n} $
(the group of permutations of a set with $ n $ elements) such that $ \psi (B_{ii})
= \mu _{i i} E_{\sigma (i) \sigma (i)}. $ Next, for the cases $ i = j $
and $ j = s, $ we have $ \sigma (i) = u(i,s) $ and $ \sigma (j) = v(i,j), $
respectively. Therefore, there exists a $ \sigma \in S_{n} $ such that,
for every pair $ i, j $ we have $ \psi (B_{ij}) =
\mu _{i j} E_{\sigma (i) \sigma (j)}. $ Now we define a $ R-$linear map
$ \rho : \ (M_{n}(R), M_{n}(R^{+})) \longrightarrow (M_{n}(R), M_{n}(R^{+})) $
by $ \rho (E_{st}) = E_{\sigma ^{-1}(s) \sigma ^{-1}(t)}. $ It is easy to
verify that $ \rho $ is an automorphism of $ l-$algebra. Let $ \tau = \rho \circ \psi , $
then $ \tau $ is an isomorphism of $ l-$algebras and $ \tau (B_{ij}) = \mu _{i j}
E_{i j}. $ By the above equality $ B_{ij} B_{js} =
q_{i j} q_{j s} q_{i s}^{-1} B_{is}, $ we get $ \tau (B_{ij}) \tau (B_{js}) =
q_{i j} q_{j s} q_{i s}^{-1} \tau (B_{is}), $ and then $ \mu _{i j} \mu _{j s}
E_{i s} = q_{i j} q_{j s} q_{i s}^{-1} \mu _{i s} E_{i s}. $ Hence
$ \mu _{i j} \mu _{j s} \mu _{i s}^{-1} = q_{i j} q_{j s} q_{i s}^{-1}, $ which gives
a positive solution for the given system of equations.

Conversely, let $ x_{ij} = \mu _{ij} \in R^{\times } $ be a positive solution
of the given system of equations. We define a $ R-$linear map
$ \varphi : \ (M_{n}(R), P) \longrightarrow
(M_{n}(R), M_{n}(R^{+})) $ by $ \varphi (B_{ij}) = \mu _{ij} E_{ij} \
(1 \leq i, j \leq n). $ Then it can be easily verified that $ \varphi $ is an
isomorphism of $ l-$algebras. The proof is completed. \quad $ \Box $
\par \vskip 0.2 cm

{\bf Lemma 2.3.} \ The system of equations $ x_{ij} x_{js} x_{is}^{-1} =
q_{ij} q_{js} q_{is}^{-1} \ (1 \leq i, j, s \leq n) $ with variables
$ x_{i^{\prime } j^{\prime }} \ (1 \leq i^{\prime }, j^{\prime } \leq n) $
in Theorem 2.2 always has positive solutions in $ R^{\times }. $
\par \vskip 0.1 cm
{\bf Proof.} \ For any $ n, x_{ii} = q_{ii} $ are always positive units of $ R $
and $ x_{ij} x_{ji} = q_{ij} q_{ji}. $ For the case $ n = 2, $ the given equations
have a solution $ x_{11} = q_{11}, x_{22} = q_{22}, x_{12} = 1, x_{21} = q_{12}q_{21}. $
We use induction on $ n \geq 3. $ If $ n = 3, $ then the given system of equations
has positive solutions in $ R^{\times } $
if and only if the following system of equations (S1) has positive solutions
in $ R^{\times }: $  \\
(S1) \ $ \left \{
   \begin{array}{l}
  x_{12}x_{21}=q_{12}q_{21}, \\
  x_{13}x_{31}=q_{13}q_{31}, \\
  x_{23}x_{32}=q_{23}q_{32}, \\
  x_{12}x_{23}x_{13}^{-1}=q_{12}q_{23}q_{13}^{-1}.
\end{array}
  \right.  $  \\
It is easy to see that (S1) has positive solutions in $ R^{\times }. $ Now for the case $ n = k, $
we assume that the given system of equations has positive solutions in $ R^{\times }. $
We want to verify the
case $ n = k + 1. $
To see this note that, just as for the case $ n = 3, $ to solve these equations
in variables $ x_{ef}, x_{fe}, x_{eg}, x_{ge}, x_{fg}, x_{gf} \ (1 \leq e < f < g \leq k + 1) $
it is sufficient to solve the following equations \\
(S2) \ $ \left \{
   \begin{array}{l}
  x_{ef}x_{fe}=q_{ef}q_{fe}, \\
  x_{eg}x_{ge}=q_{eg}q_{ge}, \\
  x_{fg}x_{gf}=q_{fg}q_{gf}, \\
  x_{ef}x_{fg}x_{eg}^{-1}=q_{ef}q_{fg}q_{eg}^{-1}.
  \end{array}
  \right.  $  \\
Using this reduction, it is easy to see that the given system of equations in case $ n = k + 1 $
has positive solutions in $ R^{\times } $ if and only if so does the
system (S3) consisting of the following four parts \\
(S3.1) \ $ x_{ij}x_{js}x_{is}^{-1} = q_{ij}q_{js}q_{is}^{-1} \
(1 \leq i, j, s \leq k), $ \\
(S3.2) \ $ \left \{
   \begin{array}{l}
x_{1, k+1} x_{k+1, 1} = q_{1, k+1} q_{k+1, 1}, \\
x_{2, k+1} x_{k+1, 2} = q_{2, k+1} q_{k+1, 2}, \\
\cdots \cdots  \\
x_{k, k+1} x_{k+1, k} = q_{k, k+1} q_{k+1, k},
\end{array}
  \right. $ \\
(S3.3) \ $ \left \{
   \begin{array}{l}
x_{12} x_{2, k+1} x_{1, k+1}^{-1} = q_{12} q_{2, k+1} q_{1, k+1}^{-1}, \\
x_{13} x_{3, k+1} x_{1, k+1}^{-1} = q_{13} q_{3, k+1} q_{1, k+1}^{-1}, \\
\cdots \cdots  \\
x_{1k}x_{k, k+1} x_{1, k+1}^{-1} = q_{1k} q_{k, k+1} q_{1, k+1}^{-1},
\end{array}
  \right. $ \\
(S3.4) \ $ x_{ij} x_{j, k+1} x_{i, k+1}^{-1} = q_{ij} q_{j, k+1} q_{i, k+1}^{-1}
\ (2 \leq i < j \leq k). $ \\
Note that the system (S3.4) can be deduced from the systems (S3.3) and (S3.1),
that is, from the equations $
x_{1j} x_{j,k+1} x_{1,k+1}^{-1} = q_{1j}q_{j, k+1}q_{1, k+1}^{-1}, \
x_{1i} x_{i, k+1} x_{1, k+1}^{-1} = q_{1i}q_{i, k+1} q_{1, k+1}^{-1} $
and $ x_{1i}x_{ij}x_{1j}^{-1}=q_{1i}q_{ij}q_{1j}^{-1}, $ we get
the equations $ x_{ij} x_{j, k+1} x_{i, k+1}^{-1} = q_{ij} q_{j, k+1} q_{i, k+1}^{-1}
\ (2 \leq i < j \leq k). $ So the given system of equations in case $ n = k + 1 $
has positive solutions in $ R^{\times } $ if and only if so does
the system consisting of (S3.1), (S3.2) and (S3.3).
By the induction hypothesis, (S3.1) has positive solutions in $ R^{\times }. $
So we can take a solution $ x_{ij} = \mu _{ij} \ (1 \leq i, j \leq k). $ So the
system consisting of (S3.1), (S3.2) and (S3.3) has a solution
$ x_{ij} = \mu _{ij} \ (1 \leq i, j \leq k), x _{1, k+1} = 1,
x _{i, k+1} = q_{1i} q_{i, k+1} q_{1, k+1}^{-1} \mu_{1, i}^{-1} \
(2 \leq i \leq k), x _{k+1, i} = x _{i, k+1}^{-1} q_{i, k+1} q _{k+1, i} \
(1 \leq i \leq k), $ which implies that the given system has positive solutions
in $ R^{\times }. $ Therefore, the conclusion holds for the case $ n = k + 1. $
The proof is completed. \quad $ \Box $
\par \vskip 0.2 cm

By the above Theorem 2.2 and Lemma 2.3, we obtain the following result,
i.e., Weinberg's conjecture holds on lattice-ordered matrix algebras over
all real GCD-domains.
\par \vskip 0.2 cm

{\bf Theorem 2.4.} \ Let $ (M_{n}(R), P) $
be an $ l-$algebra over $ R $ with $ I \in P, $ then
$ (M_{n}(R), P) \cong (M_{n}(R), M_{n}(R^{+})). $
\par \vskip 0.2 cm

{\bf Theorem 2.5.} \ Any $ l-$algebra $ (M_{2}(R), P) $
is isomorphic to $ (M_{2}(R), P_{A}), $ where
$ P_{A} = A M_{2}(R^{+}) $ for some $ A \in M_{2}(R^{+}) $
with $ \text{det}(A) \in R^{\times}. $
\par \vskip 0.1 cm
{\bf Proof.} \ By Proposition 2.1 above, $ (M_{2}(R), P) $
has a $ vl-$basis $ B= \{B_{ij} \} \ (1 \leq i,j \leq 2) $ with
$ B_{ij} = q_{ij} H D E_{ij} H^{-1} $ and $ B_{ij} B_{rs} = d_{jr} q_{ij}
q_{rs} q_{is}^{-1} B_{is}. $ From [MW],
$ (M_{2}(K), \overline{P}) \cong (M_{2}(K), P_{D}) $
where $ P_{D} = D M_{2}(K^{+}) $ for one of the following
three matrices $ D: $ \\
(1) \ $ D = I; \quad (2) \  D =  \begin{pmatrix}  1 & 1 \\
1 & 0
\end{pmatrix};  \quad
(3) \ D = \begin{pmatrix}
1 & 1 \\
a & b
\end{pmatrix}, $ where $ a, b \in K $ and $ a > b > 0. $ \\
First, for the case (1), the conclusion follows directly from the above
Theorem 2.4.

Next for the case (2), by Proposition 2.1.(2), we have
$ d_{jr} q_{ij} q_{rs} q_{is}^{-1} \in R^{+} \
(1 \leq i, j, r, s \leq 2). $ From $ (i, j , r, s) = (1, 1, 1, 1),
(1, 2, 1, 1), (1, 1, 2, 1), (1, 1, 2, 2), $ we get
$ q_{11}, q_{12}, q_{21}, q_{11}q_{22} \in R^{+}, $ respectively.
By Proposition 2.1.(4), we have $ q_{11} q_{12} q_{21} q_{22} \in R^{\times }. $
So $ q_{12}, q_{21}, q_{11}q_{22} \in R^{\times } \cap R^{+}. $
Set $ A = \begin{pmatrix}
  q_{11} & q_{21} \\
  q_{12} & 0
\end{pmatrix} $ and $ x_{11} = x_{12} = x_{21} = 1, \ x_{22} =
q_{11} q_{22} q_{12}^{-1} q_{21}^{-1}. $ Obviously, $ det(A) \in R^{\times}. $
Also it is well known that, for $ B \in M_{n}(R^{+}),  B M_{n}(R^{+}) $
is the positive cone of a lattice order on $ M_{n}(R) $
if and only if $ \text{det}(B) \in R^{\times } $ (see [St2], p.595).
So we know that $ (M_{2}(R), P_{A}) $
is an $ l-$algebra with positive cone $ P_{A} = A M_{2} (R^{+}). $
We define a $ R-$linear map
$ \phi : \ (M_{2}(R), P) \longrightarrow
(M_{2}(R), P_{A}) $ by $ \phi (B_{ij}) = x_{ij} A E_{ij} \
(1 \leq i, j \leq 2). $ By direct computation of various subcases that
arise from the four cases $ (j, r) = (2, 2), (1, 1), (1, 2), (2, 1), $ it can be seen
that $ \phi $ preserves multiplication. So $ \phi $ is an isomorphism of
$ l-$algebras.

For the last case (3), just as done in the case (2) above, by Proposition 2.1 (2) (4),
we have $ a q_{12}, b q_{22}, q_{11}, q_{21} \in R^{+}. $
We assert that $ q_{11} q_{22} = \mu q_{12} q_{21} $ for some
$ \mu \in R^{+} \cap R^{\times } $ (we will prove it later).  Set
$ C = \begin{pmatrix}
  \frac{q_{11}}{p_{11}} & \frac{q_{21}}{p_{21}} \\
  \frac{aq_{12}}{p_{12}} & \frac{bq_{22}}{p_{22}}
\end{pmatrix}, $ where $ p_{ij} \in R^{\times } \cap R^{+} $ satisfy
$ \frac{q_{11} q_{22}}{p_{11} p_{22}} = \frac{q_{12} q_{21}}{p_{12} p_{21}}.
$ Let $ t = \frac{q_{12} q_{21}}{p_{12} p_{21}}, $ then $ C \in M_{2}(R^{+})
$ with $ \text{det}(C) = (b-a) t \in R. $ By Proposition 2.1, $
(q_{12} q_{22}q_{11}q_{21})(\text{det}(D))^{2} = (q_{12} q_{22} q_{11} q_{21})
(b-a)^{2} \in R^{\times}, $ so $ t^{2}(b-a)^{2} \in R^{\times }, $ hence
$ \text{det}(C) \in R^{\times }. $ Then
$ (M_{2}(R), P_{C}) $ is an $ l-$algebra with positive cone
$ P_{C} = C M_{2} (R^{+}). $ We define a $ R-$linear map
$ \theta : \ (M_{2}(R), P) \longrightarrow
(M_{2}(R), P_{C}) $ by $ \theta (B_{ij}) = p_{ij} C E_{ij} \
(1 \leq i, j \leq 2). $ Also by direct computation of various subcases, it can be seen
that $ \theta $ preserves multiplication. So $ \theta $ is an isomorphism of
$ l-$algebras. Now we only need to prove our assertion
$ q_{11} q_{22} = \mu q_{12} q_{21} $ for some $ \mu \in R^{+} \cap R^{\times }. $
Note that $ a q_{12}, b q_{22}, q_{11}, q_{21} \in R^{+} $ and
$ D = \begin{pmatrix}
  1 & 1 \\
  a & b
\end{pmatrix}. $ Let $ I = \sum_{i, j=1}^{2} k_{i j} B_{ij}, \ k_{ij} \in R, $
then $ I = k_{11} q_{11} D E_{11} + k_{12} q_{12} D E_{12} + k_{21} q_{21} DE_{21}
+ k_{22} q_{22} D E_{11}. $ By direct calculation, $ I = - \frac{b}{a-b} D E_{11} +
\frac{1}{a - b} D E_{12} + \frac{a}{a - b} D E_{21} - \frac{1}{a - b} D E_{22}. $
So $ k_{11} q_{11} = - \frac{b}{a - b}, \ k_{12} q_{12} = \frac{1}{a-b}, \
k_{21} q_{21} = \frac{a}{a - b}, \ k_{22} q_{22} = - \frac{1}{a - b}. $
Let $ m = \frac{a}{a - b}, \ m - 1 = \frac{b}{a - b}, $ then $
k_{21} q_{21} = a k_{12} q_{12} = m, \ k_{11} q_{11} = b k_{22} q_{22} = 1 - m. $
Let $ \varepsilon = (q_{12} q_{22} q_{11} q_{21})(\text{det}(D))^{2}, $ then
by Proposition 2.1, $ \varepsilon \in R^{\times }. $ Since $ R $ is a GCD-domain and
$ \text{gcd}(m, m-1) = 1, $ we have $$
(\varpi ): \ \text{gcd}(m, q_{11}) = 1, \ \text{gcd}(m, bq_{22}) = 1, \
\text{gcd}(m-1, q_{21}) = 1, \ \text{gcd}(m-1, a q_{12}) = 1. $$
Since $ q_{11} a q_{12} q_{21} b q_{22} = \frac{a b}{(a - b)^{2}}q_{11} q_{12}q_{21}
q_{22} (\text{det}(D))^{2} = \frac{a b \varepsilon }{(a-b)^{2}}
=(m - 1) m \varepsilon , $ by the equalities $(\varpi ), $ we get
$ b q_{11} q_{22} = (m - 1) \mu_{1}, \ a q_{12} q_{21} = m \mu_{2}, $
where $ \mu_{1}, \mu_{2} \in R^{\times} $ and $ \mu_{1} \mu_{2} = \varepsilon . $
Since $$ \frac{a q_{12} q_{21}}{b q_{11} q_{22}} = \frac{m \mu_{2}}{(m - 1) \mu_{1}}
\Longrightarrow \frac{a}{b} \frac{q_{12} q_{21}}{q_{11} q_{22}} =
\frac{a}{a - b} \frac{a - b}{b} \frac{\mu_{2}}{\mu_{1}}, $$
we obtain that $ \mu = \frac{q_{12} q_{21}}{q_{11} q_{22}} = \frac{\mu_{2}}{\mu_{1}}, $
i.e., $ q_{11} q_{22} = \mu q_{12} q_{21}, \ \mu\in R^{\times}. $
This proves our assertion, and the proof is completed. \quad $ \Box $

\par  \vskip 0.3 cm

{ \bf Acknowledgments. } \ We are grateful to Professor
R. H. Redfield for sending us helpful materials. We thank
Professor Jingjing Ma for helpful suggestions. We also thank the anonymous
referee for helpful suggestions, especially for suggesting us to replace
UFD rings by GCD-domains in this paper.

\par  \vskip 0.2 cm

\hspace{-0.8cm} {\bf References }

\begin{description}

\item[[Bi]] G. Birkhoff, Lattice theory, 3rd Edition,
Amer. Math. Soc. Colloquium Pub. 25, Amer. Math. Soc.,
Providence, 1995.
\item[[K]] I. Kaplansky, Commutative Rings, Polygonal Publishing House,
Washington, New Jersey, 1994.
\item[[MR]] J. Ma, R. H. Redfield, Lattice-ordered matrix
rings over the integers, Comm. Algeba. 35 (2007), 2160-2170.
\item[[MW]] J. Ma, P. Wojciechowski, Lattice orders on matrix
algebras, Algebra Universalis, 47(2002), 435-441.
\item[[St1]] S.A.Steinberg, On the scarcity of lattice-ordered matrix
algebras II, Proc. Amer. Math. Soc., 128 (2000), 1605-1612.
\item[[St2]] S.A.Steinberg, Lattice-ordered Rings and Modules, Springer,
New York, 2010.
\item[[W]] E. Weinberg, On the scarcity of lattice-ordered matrix rings,
Pacific J. Math. 19 (1966), 561-571.

\end{description}

\end{document}